\documentclass[a4paper]{amsart}

\usepackage{graphics,amssymb,enumerate}
\usepackage{hyperref}
\usepackage[dvips]{graphicx}
\usepackage[all]{xy}

\newtheorem{theorem}{Theorem}

\newtheorem{cor}[theorem]{Corollary}

\newtheorem{rem}[theorem]{Remark}
\def\bd{\begin{displaymath}}
\def\ed{\end{displaymath}}
\def\be{\begin{equation}}
\def\ee{\end{equation}}
\def\p{\partial}

\def\g{\mathfrak{g}}
\def\h{\mathfrak{h}}
\newcommand{\mc}{\mathcal}
\def\pb#1#2{\left\{#1,#2\right\}}
\def\Pb{\pb{\cdot}{\cdot}}

\def\R{\mbox{$\mathbb R$}}

\def\set#1{\left\{#1\right\}}
\newcommand{\tr}{\text{\rm tr}\,}

\begin{document}

\title[Toda  Volterra]{Reduction and Realization in Toda and Volterra}

\author{Pantelis. A. ~Damianou}
\address{Department of Mathematics and Statistics\\
University of Cyprus\\
P.O.~Box 20537, 1678 Nicosia\\Cyprus} \email{damianou@ucy.ac.cy}

\thanks{Dedicated to H. Flaschka and J. Moser}

\begin{abstract}
We construct a new symplectic,  bi-hamiltonian realization of the KM-system by
reducing the corresponding one for the  Toda lattice. The bi-hamiltonian pair is  constructed
using a reduction theorem of Fernandes and Vanhaecke. In this paper we also review  the important work of Moser on the Toda and KM-systems.
\end{abstract}

\date{}
%%% ----------------------------------------------------------------------
\maketitle
%%% ----------------------------------------------------------------------

\section{Introduction} \label{intro}

Some of the most important work of Moser concerns the Toda and Volterra lattices. These two systems are closely related and Moser gave an explicit construction demonstrating the relationship between the two systems. The Toda lattice was discovered by Morikasu Toda in 1967. Computer experiments by Ford et al. \cite{ford} suggested that the Toda lattice is integrable. In two papers, both in the same issue of Physical Review B (\cite{flaschka1}, \cite{henon})  H\'enon and Flaschka in 1974 demonstrated the integrability of the system. H\'enon provided   the required independent constants of motion using combinatorial arguments.  Flaschka proved integrability of the  lattice via a change of variables and the construction of a Lax pair.  A few months later, Manakov \cite{manakov}  established independently similar results.
The Toda lattice is a discrete approximation of the KdV equation. In \cite{flaschka2} Flaschka used this analogy to solve the system  using a discretization of the inverse scattering method of Gardner, Greene, Kruskal and Miura. The method which was developed for  partial differential equations was applied  for the first time for  the infinite Toda lattice  with success.
Moser's paper \cite{moser1} deals with the finite non-periodic Toda lattice. Moser solves the system
by constructing action-angle coordinates.  His action variables are the eigenvalues of the Lax matrix  of Flaschka and the angles are the last components of the  normalized eigenvectors. We will describe the solution of the Toda lattice by Moser in Section 2.
The Volterra system  (also known as the KM system) is closely related to the Toda lattice, and the connection will become even more intimate after the results of the present paper. It was first solved by
Kac and van-Moerbeke in \cite{kac}, using the  discrete version of
inverse scattering of Flaschka.  Moser in \cite{moser2}  solved this system explicitly  in a fashion similar to \cite{moser1}, i.e.  by construction of action-angle coordinates, using the Weyl function and continued fractions. We will not examine this result in the paper since it is similar (although more complicated) to the corresponding one for the Toda lattice which we describe in Section 2.

The purpose of this paper is to understand completely the following diagram:

\bd
\xymatrix{
{\rm Toda} (q,p) \ar[rr]^{F}\ar[dd]_{Invol.}
&&
{\rm Toda} (a, b) \ar[dd]^{Invol.}\\
\\
{\rm Volterra} (q) \ar[rr]^{G}&& {\rm Volterra} (a)
}
\ed

The top arrow is the well-known symplectic realization of the Toda lattice where $F$ is the Flaschka transformation.   We can think of it as the following symplectic bi-hamiltonian realization:

\bd
F: (J_1, J_2, h_1, h_2)  \to  (\pi_1, \pi_2, H_1, H_2) \ ,
\ed
 where $\pi_1$ and $\pi_2$ are the well-known linear and quadratic Toda brackets and $H_1$, $H_2$ are the trace of the Lax matrix and the Hamiltonian respectively. They form a bi-Hamiltonian pair: $\pi_1 dH_2= \pi_2 dH_1$.  Similarly, $J_1$ is the standard symplectic bracket while $J_2$ is the Das-Okubo bracket in $(q, p)$ coordinates. $h_1$ and $h_2$ correspond to the sum of the momenta and the Hamiltonian respectively. Here, we also have a bi-hamiltonian formulation
 $J_1 dh_2=J_2 dh_1$.  The multi-hamiltonian structure in $(a,b)$  coordinates was obtained  in \cite{damianou1}, \cite{damianou3} while in $(q,p)$  coordinates in \cite{fern1} using a Theorem of Oevel. There is a recursion operator defined by ${\mc R}=J_2 J_1^{-1}$.

 The bottom arrow is a new symplectic  bi-hamiltonian realization of the KM-system.  Let us denote it by
 \bd
 G: (w_2, w_3, i_1, i_0) \to (v_2, v_3, I_1, I_0) \ .
 \ed

 The brackets $v_2$ and $v_3$ are the well-known quadratic and cubic brackets for the KM-system defined in \cite{damianou2}. The functions $I_i$ are defined by  $I_1=\frac{1}{2} {\rm tr}L^{2}$, and $I_0=\log |det(L)|$,  where $L$ is the Lax matrix. The phase space for the KM-system $V_{a}$ is of odd dimension  $(a_1, a_2, \dots, a_{2n+1})$. Letting $N=2(n+1)$ the dimension of $V_a$ is $N-1$, while the dimension of $T_{ (a, b)}$ is $2N-1$ and the  dimension of $T_{(q,p)}$  is $2N$.

  We consider the   space $V_{q}$  in the  variables $(q_1, q_2, \dots, q_N)$.
  On this space we  define  a symplectic Poisson bracket $w_2$ by the formula

 \bd
 \{ q_i, q_j \}=1  \qquad \forall \ \  i < j  \ .
 \ed
 This is a constant symplectic bracket and the matrix $\{q_i, q_j \}$ has  determinant one.

This
bracket $w_2$ corresponds to the quadratic
bracket $v_2$ via  the mapping
\bd G(q_1, \dots, q_N) =\left( e^{q_1- q_2}, \dots, e^{q_{N-1}- q_N} \right) \ . \ed

   We then
define a bracket $w_3$ in ${\bf R}^{N}$ which is mapped to
$v_3$ under the transformation $G$.  The  bracket $w_3$ is defined by the formula:

\bd \{ q_i, q_j \}=e^{q_{i-1}-q_i}+(1-\delta_{i+1, j}) e^{q_i- q_{i+1}} + e^{q_{j-1}-q_j}+e^{q_j-q_{j+1}} \ . \ed
Whenever a term is not defined we omit that whole term.
Define
\bd i_0=q_1-q_2 + q_3-q_4+ \dots +q_{N-1}-q_N  \   \ed
and the  Hamiltonian $i_1$  given by

\bd
i_1=  \sum_{i=1}^{N-1}
e^{q_{i}-q_{i+1}}  \ .
  \ed

Then we have the following  bi-hamiltonian pair:

\bd w_2 di_1= w_3 di_0 \ . \ed

 Finally, let us comment the two vertical arrows of the diagram. The right hand arrow is described in \cite{damianou5} where the multi-hamiltonian structure of the KM-system is obtained by reducing the Toda hierarchy using a Theorem of Fernandes and Vanhaecke \cite{fern3}.   The left hand arrow is similarly a reduction of the multi-hamiltonian structure of the Toda lattice in $(q,p)$ coordinates  to a new multi-hamiltonian structure for the KM-system in $V_q$ space.  We use the involution $\psi (q, p)=(q,-p)$ which is a  Poisson automorphism of both $J_2$ and $J_4=\mc{R}^2 J_2$.  The reduction of  $J_2$ and $J_4$ produces the brackets  $w_2$ and $w_3$ respectively.  We should mention that Moser had an algorithm for going in the opposite direction of the right hand side arrow, i.e.   a procedure for obtaining Toda equations from the Volterra equations by squaring the Lax matrix and some chopping.  This is described briefly in Section 10.

\section{Moser's solution of the Toda lattice}

The Toda lattice  with Hamiltonian function

\be \label
{toda}
  H(q_1, \dots, q_N, \,  p_1, \dots, p_N) = \sum_{i=1}^N \,  { 1 \over 2} \, p_i^2 +
\sum _{i=1}^{N-1} \,  e^{ q_i-q_{i+1}}  \ ,
\ee

can be transformed via  a change of variables due to Flaschka \cite{flaschka1}  to a Lax  pair of the form   $\dot L = [B, L] $, where $L$ is the Jacobi matrix

\be
 L= \begin{pmatrix} \label{jacobi} b_1 &  a_1 & 0 & \cdots & \cdots & 0 \cr
                   a_1 & b_2 & a_2 & \cdots &    & \vdots \cr
                   0 & a_2 & b_3 & \ddots &  &  \cr
                   \vdots & & \ddots & \ddots & & \vdots \cr
                   \vdots & & & \ddots & \ddots & a_{N-1} \cr
                   0 & \cdots & & \cdots & a_{N-1} & b_N     \end{pmatrix} \ ,
\ee

and $B$ is the projection of $L$ into the skew-symmetric part of $L$ in the Lie algebra decomposition, lower triangular plus skew-symmetric.
This is an example of an isospectral deformation; the entries of $L$ vary over time but the eigenvalues  remain constant. It follows that the  symmetric polynomials of the eigenvalues,  $ H_i={1 \over i} {\rm tr} \, L^i$ are  constants of
motion.
The differential equations in the case of symmetric tri-diagonal Lax matrix are:
\be
\begin{array}{lcl}
 \dot a _i& = & a_i \,  (b_{i+1} -b_i )    \\
   \dot b _i &= & 2 \, ( a_i^2 - a_{i-1}^2 ) \ .  \label{toda-tri}
\end{array}
\ee

Moser's  elegant  solution  of the open Toda lattice uses the
Weyl function $f(\lambda)$ and an old (19th century) method of
Stieltjes which connects the continued fraction  of $f(\lambda)$
with its partial fraction expansion.
 The key ingredient is the
map which takes the $(a,b)$ phase space of tridiagonal Jacobi
matrices to a new space of variables $(\lambda_i, r_i)$
 where $\lambda_i$ is an eigenvalue of the Jacobi matrix and $r_i$ is  the residue of $f(\lambda)$.
 We present a brief outline of Moser's construction.

Moser in \cite{moser1} introduced the resolvent \bd
R(\lambda)=\left( \lambda I - L \right)^{-1} \ , \ed

 and defined the Weyl function
\bd f(\lambda)=R_{NN}(\lambda)=\left( R(\lambda)e_N, \  e_N \right) \ , \ed

where $e_N=(0,0, \dots, 0, 1)$.

The function $f(\lambda)$ has a simple pole at $\lambda=
\lambda_i$ and so it admits a partial fraction expansion

\begin{equation} \label{pf}
f(\lambda)= \sum_{i=1}^{N} \frac{r_i^2}{\lambda-\lambda_i}
  \ ,
\end{equation}

with positive residue $r_i^2$.   It is clear from the formula for calculating the inverse of a matrix that
$f(\lambda)=\frac{\Delta_{N-1} } {\Delta_N}$.  We denote by $\Delta_k$ the $k$ by $k$ sub-determinant obtained by deleting the last $N-k$ rows and columns of $\lambda I- L$.  It follows that $\lim_{\lambda \to \infty } \lambda f(\lambda)=1$ and therefore the residue of $f(\lambda)$ at infinity is $-1$.  As a result,
\bd \sum_{i=1}^N r_i^2=1 \ . \ed
More generally, one has a recursion formula:
\be \label{recursion-formula}   \Delta_k=(\lambda- b_k) \Delta_{k-1} -a_{k-1}^2 \Delta_{k-2} \ . \ee
Moser notes that the mapping $\phi$ between

\bd T_{(a,b)}=\{ (a_1, \dots, a_{N-1}, b_1, \dots, b_N,  \ \ \  a_i > 0 ) \} \ed

and
\bd T_{(\lambda, r)}  =\{ (\lambda_1, \dots, \lambda_N, r_1, \dots, r_N , \ \ \lambda_1 < \lambda_2 < \dots < \lambda_N,  \ \ \sum_{i=1}^N r_i^2 =1, \ \  r_i > 0) \} \ed
is one  to one and onto.  The inverse mapping $\phi^{-1}\ : \  T_{(\lambda, r)}  \to T_{(a,b)}  $ corresponds to the inverse scattering transform.

 Moser derives the  differential equations in the variables $(\lambda, r)$
using the Lax pair of the system and the recursion formula (\ref{recursion-formula}).  The equations take the form

\begin{eqnarray}
\dot{\lambda}_i& = & 0 \nonumber\\
\dot{r}_i & = & -(\lambda_i- \sum_{j=1}^N \lambda_j r_j^2)   r_i\;. \label{a7}
\end{eqnarray}

If one considers  $r_i$ as homogeneous variables then the differential equations (\ref{a7})  become linear
\bd  \dot{\lambda}_i =  0  \qquad  \dot{r_i}=-\lambda_i r_i \ . \ed

 The variables $a_i^2$, $b_i$ may be expressed as rational
functions of $\lambda_i$ and $r_i$ using a continued fraction
expansion of $f(\lambda)$ which dates back to Stieltjes. Since the
computation of the continued fraction from the partial fraction
expansion is a rational process the solution is expressed as a
rational function of the variables $(\lambda_i, \ r_i)$. The
procedure is as follows:

The $R_{NN}$ element of the resolvent, as defined previously,
takes the following continued fraction representation:

\begin{equation}
f(\lambda)=\frac{1}{\lambda-b_N-\frac{a_{N-1}^2}{\lambda-b_{N-1}-\frac{a_{N-2}^2}{\frac{\vdots}
{\lambda-b_2-\frac{a_1^2}{\lambda-b_1}}}}}   \ .
\end{equation}

 Stieltjes described a procedure that
allows one to  express $a_i$ and $b_i$ in terms of $\lambda_1,
\ldots,\lambda_N$ and $r_1,\ldots,r_N$. We briefly describe the
method. We expand the partial fraction expansion of $f(\lambda)$
 given in (\ref{pf}) in a series of powers of
$\frac{1}{\lambda}$.

 The coefficient  of $\lambda^{j+1}$ is denoted by $c_j$ and
equals,

\bd c_j= \sum_{i=1}^{N} r_i^2 \lambda_i^j
, \quad j=0, 1, \ldots. \ed The formulas of
Stieltjes involve certain $i \times i$ determinants which we now
define:

\bd
A_i=\left| \begin{array}{cccc} c_0 & c_1 & \ldots & c_{i-1} \\
                    c_1 & c_2 & \ldots & c_{i} \\
                    \vdots & & &\\
                    c_{i-1} & c_{i} & \ldots & c_{2i-2}
    \end{array} \right|, \qquad B_i=\left| \begin{array}{cccc} c_1 & c_2 & \ldots & c_{i} \\
                    c_2 & c_3 & \ldots & c_{i+1} \\
                    \vdots \\
                    c_{i} & c_{i+1} & \ldots & c_{2i-1}
    \end{array} \right|.
\ed The formulas that give the relation between the variables
$(a,b)$ and $(r,\lambda)$ are,

\begin{eqnarray*}
&& a_{N-i}^2=\frac{A_{i-1} A_{i+1}}{A_i^2}, \hspace{2.32cm}
i=1,\ldots,N-1
\\
&&  b_{N+1-i}=\frac{A_i
B_{i-2}}{A_{i-1}B_{i-1}}+\frac{A_{i-1}B_{i}}{A_i B_{i-1}}, \quad
i=1,\ldots,N
\end{eqnarray*}
where $A_0=1$, $B_0=1$, $B_{-1}=0$.

For example, in the case  $N=2$

\bd A_1=c_0, \quad A_2=c_0 c_2-c_1^2, \quad B_1=c_1, \quad B_2=c_1
c_3-c_2^2 \ed and therefore \bd a_1^2=A_2, \quad
b_1=\frac{A_2}{B_1}+\frac{B_2}{A_2 B_1}, \quad b_2=B_1 \ . \ed
Thus,
\begin{eqnarray}
&& a_1^2=\frac{r_1^2 r_2^2 (\lambda_2-\lambda_1)^2}
{(r_1^2+r_2^2)^2}\nonumber
\\
&& b_1=\frac{r_1^2 \lambda_2+r_2^2 \lambda_1}{r_1^2+r_2^2}  \label{t1}\\
&& b_2=\frac{r_1^2 \lambda_1+r_2^2 \lambda_2}{r_1^2+r_2^2} \
.\nonumber
\end{eqnarray}
We have written the solution in this form to show that $r_i$ are homogeneous coordinates, i.e.  the solution does not change if we replace $r_i$ by  a non-zero multiple.   As we mentioned earlier in these homogeneous coordinates the differential equations become linear.
One can check that the differential equations (\ref{a7})  correspond via transformation (\ref{t1})
to the $A_2$ Toda equations

\begin{eqnarray*}
&& \dot{a}_1=a_1(b_2-b_1)  \\
&& \dot{b}_1=2 a_1^2   \\
&& \dot{b}_2=-2 a_1^2 \ .
\end{eqnarray*}

As Moser notes, it is not too hard to obtain explicit expressions
for $N=3$ but the general case is quite complicated.

\section{Background}

 Consider a
differential equation on a manifold $M$  defined by a vector field
$\chi$.  A  vector field $Z$ is a   symmetry of the equation  if
\begin{equation*}
[Z, \chi]=0  .
\end{equation*}
A vector field $Z$ is called a master symmetry if
\begin{equation*}
[[Z, \chi], \chi]=0 ,
\end{equation*}
but
\begin{equation*}
[Z, \chi] \not= 0  .
\end{equation*}

\noindent Master  symmetries were first introduced by Fokas and
Fuchssteiner in \cite{fokas1}  in connection with the Benjamin-Ono
Equation.

A bi-Hamiltonian system is defined by specifying two Hamiltonian
functions $H_1$, $H_2$ and two Poisson tensors $\pi_1$ and
$\pi_2$, that give rise to the same Hamiltonian equations. Namely,
$
%\begin{equation*}
\pi_1 \nabla H_2=\pi_2 \nabla H_1.
%\end{equation*}
$ The notion of bi-Hamiltonian structures is due to Magri
\cite{magri}. Suppose that we have a bi-Hamiltonian system defined
by the Poisson tensors $\pi_1$, $\pi_2$ and the Hamiltonians
$H_1$, $H_2$.
 Assume that $\pi_1$ is symplectic.  We define
the recursion operator ${\mathcal R} = \pi_2 \pi_1^{-1}$,  the
higher flows
\begin{equation*}
\chi_{i} = {\mathcal R}^{i-1} \chi_1 \ ,
\end{equation*}
and the higher order Poisson tensors \bd \pi_i = {\mathcal
R}^{i-1} \pi_1 \ . \ed

\noindent For a non-degenerate bi-Hamiltonian system, master
symmetries can be generated using a method due to Oevel
\cite{oevel}.

\begin{theorem} \label{oevel}
Suppose that   $X_0$ is a conformal symmetry for both  $\pi_1$,
$\pi_2$ and $H_1$, i.e.  for some scalars $\lambda$, $\mu$, and
$\nu$ we have \bd {\mathcal L}_{X_0} \pi_1= \lambda \pi_1,
\quad{\mathcal L}_{X_0} \pi_2 = \mu \pi_2, \quad {\mathcal
L}_{X_0} H_1 = \nu H_1  . \ed Then the vector fields $X_i =
{\mathcal R}^i X_0$ are master symmetries and we have,
\bd
\begin{array}{lcl}
 (a) \ {\mathcal L}_{X_i} H_j & =& (\nu +(j-1+i) (\mu -\lambda))
H_{i+j}
\cr
 (b) \ {\mathcal L}_{X_i} \pi_j &= &(\mu +(j-i-2) (\mu -\lambda))
\pi_{i+j} \cr
 (c) \ [X_i, X_j] &=& (\mu - \lambda) (j-i) X_{i+j}  \ .
\end{array}
\ed

\end{theorem}

\smallskip

A symplectic realization (see \cite{weinstein}) of a Poisson manifold $(M, \pi)$ is a symplectic manifold $(S, \omega)$ together with a surjective Poisson submersion $f: S \to M$.  In this paper our realizations will have additional structure; they will be symplectic bi-hamiltonian realizations as in \cite{panasyuk}, \cite{petalidou}.   Suppose on $M$ we have a bi-hamiltonian pair $\pi_1$, $\pi_2$ and two  Hamiltonians $H_1$, $H_2$ such that $\pi_1 dH_1= \pi_2 dH_2$.  A symplectic bi-hamiltonian realization will be a manifold $S$, of even dimension,  together  with two Poisson tensors $J_1$, $J_2$ with $J_1$ symplectic and a surjective  submersion $F: S \to M$  which is a Poisson mapping between $J_i$ and $\pi_i$, i.e. $\{F^*f, F^*g \}_{J_i}=F^*\{f, g\}_{\pi_i}$ for all $f, g \in C^{\infty}(M)$.  In addition, $J_1 dh_1= J_2 dh_2$ where $h_i=H_i \circ F$.

We shall see that the relation between the Toda and Volterra systems
relies on special symmetries of the phase spaces.
We will need a theorem of Fernandes and Vanhaecke which gives
 conditions under which the fixed point set of a Poisson action
inherits a Poisson bracket.

Although we will be interested mainly in finite symmetries, we have
the following general result \cite{fern3}:

\begin{theorem}
\label{thm:invol:poisson}
  Suppose that $(M,\Pb)$ is a Poisson manifold, and $G$ is a compact
  group acting on $M$ by Poisson automorphisms. Let $N=M^G$ be the
  submanifold of $M$ consisting of the fixed points of the action and
  let $\iota:N\hookrightarrow M$ be the inclusion. Then $N$ carries a
  (unique) Poisson structure $\Pb_N$ such that
  \begin{equation}
    \label{eq:invol:poisson}
    \imath^*\pb{F_1}{F_2}=\pb{\imath^*F_1}{\imath^*F_2}_N
  \end{equation}
  for all $G$-invariant functions $F_1,F_2\in C^\infty(M)$.
\end{theorem}

\begin{rem}
The previous result can be seen as a particular case of Dirac
reduction (for the general theorem on Dirac reduction, see Weinstein
(\cite{weinstein}, Prop.~1.4) and Courant (\cite{courant}, Thm.~3.2.1).
\end{rem}

Since this result applies in particular when $G$ is a finite group, we have:

\begin{cor}
  Suppose that $(M,\Pb)$ is a Poisson manifold, and $G$ is a finite group
  acting on $M$ by Poisson automorphisms. Then the fixed point set $N=M^G$
  carries a (unique) Poisson structure $\Pb_N$ satisfying 
  (\ref{eq:invol:poisson}).
\end{cor}

Let us consider the special case $G={\bf Z}_2$. Then $G=\set{I,\phi}$,
where $\phi:M\to M$ is a Poisson involution. We conclude that
$N=M^G=\set{x:\phi(x)=x}$ has a unique Poisson bracket satisfying
equation (\ref{eq:invol:poisson}).  So we see that Theorem
\ref{thm:invol:poisson} contains as a special case the following
result, which is known as the \emph{Poisson involution theorem} (see
\cite{fern3, xu}).

\begin{cor}
\label{cor:Poisson:involution}
  Suppose that $(M,\Pb)$ is a Poisson manifold, and $\phi:M\to M$ is a
  Poisson involution. Then the fixed point set $N=\set{x\in M:\phi(x)=x}$
  carries a (unique) Poisson structure $\Pb_N$ such that
  \[
    \imath^*\pb{F_1}{F_2}=\pb{\imath^*F_1}{\imath^*F_2}_N
  \]
  for all functions $F_1,F_2\in C^\infty(M)$ invariant under $\phi$.
\end{cor}

\section {KM-system}

The Volterra system, also known as KM system is defined by

\begin{equation}
\dot a_i = a_i(a_{i+1}-a_{i-1}) \qquad i=1,2, \dots,n, \label{a1}
\end{equation}
where $a_0 \!= a_{n+1} \!=0$. It was studied originally by
Volterra in \cite{volterra} to describe population evolution in a
hierarchical system of competing species.  In
\cite{moser2} Moser gave a solution of the system using the method
of continued fractions and in the process he constructed
action-angle coordinates. Equations (\ref{a1}) can be considered
as a finite-dimensional approximation of the Korteweg-de Vries
(KdV) equation. They also appear in the discretization of
conformal field theory; the Poisson bracket  for this system can
be thought as a lattice generalization of the Virasoro algebra
\cite{fadeev2}. The variables $a_i$  are an intermediate step in
the construction of the action-angle variables for the Liouville
model on the lattice.

The Volterra system is usually
associated with a simple Lie algebra of type $A_n$. Bogoyavlensky
generalized this system for each simple Lie algebra and showed
that the corresponding systems are also integrable. See
\cite{bog1,bog2} for more details. The relation between Volterra
and Toda systems is also examined in \cite{damianou5}, \cite{damianou6}.

The Hamiltonian description of system (\ref{a1}) can be found in
\cite{fadeev} and \cite{damianou2}.  The Lax
pair is given by

\begin{equation*}
\dot{L}=[B, L],
\end{equation*}
where

\be \label{lax-volterra}
L= \begin{pmatrix} 0 &  1 & 0 & \cdots & \cdots & 0 \cr
                   a_1 & 0 & 1 & \ddots &    & \vdots \cr
                   0 & a_2 & 0 & \ddots &  &  \vdots \cr
                   \vdots & \ddots & \ddots & \ddots & & 0 \cr
                   \vdots & & & \ddots & \ddots & 1 \cr
                   0 & \cdots &  \cdots & 0 & a_{n} & 0  \end{pmatrix}  \ .
\ee

and

\begin{equation*}
 B=\begin{pmatrix} 0 &  1 & 0 & \cdots & \cdots & 0 \cr
                   0 & 0 & 1 & \ddots &    & \vdots \cr
                   -a_1 a_2 & 0 & 0 & \ddots &  &  \vdots \cr
                   \vdots & -a_2 a_3 & \ddots & \ddots & & 0 \cr
                   \vdots & & & \ddots & \ddots & 1 \cr
                   0 & \cdots &  \cdots & -a_{n-1} a_n & 0 & 0  \end{pmatrix} \ .
 \end{equation*}
 It
follows  that the  functions $ I_i={1 \over 2i} {\rm Tr} \, L^{2i}$ are
constants of motion. The system (\ref{a1}) is integrable for any value of $n$ but in this paper we restrict our attention only to the case where $n$ is odd.

Following \cite{damianou2} we define the following quadratic
Poisson bracket, \bd \{a_i, a_{i+1} \}=a_i a_{i+1}, \ed and all
other brackets equal to zero. We denote this bracket by $v_2$.
For this bracket det$L$ is a Casimir and the eigenvalues of $L$
are in involution. Of course, the functions $I_i$ are also in
involution. Taking the function $I_1=\sum_i^{n} a_i $ as the
Hamiltonian we obtain equations (\ref{a1}),  i.e. $v_2 dI_1$ is the KM Hamiltonian vector field.

In \cite{damianou2} one also finds a cubic Poisson bracket which
corresponds to the second KdV bracket in the continuum limit. It
is defined by the formulas,

\begin{equation*}
\begin{array}{rcl}
\{ a_i, a_{i+1} \}&=& a_i a_{i+1} (a_i+ a_{i+1}) \\
\{ a_i, a_{i+2} \} &=& a_i a_{i+1} a_{i+2} \, ,
\end{array}
\end{equation*} all other brackets are zero. We denote this bracket by
$v_3$.  In this bracket we still have involution of invariants.
We also have Lenard type relations of the form \bd v_3\, d I_{2i} = v_2 \, d I_{2i+2} . \ed

In \cite{damianou2} appears a  bracket that is homogeneous of
degree one, a rational bracket constructed using a master
symmetry. This bracket, denoted by $v_1$, has $I_1=\frac{1}{2} {\rm Tr} L^2$ as
Casimir and the Hamiltonian is $I_2={ 1\over 4} {\rm Tr} L^4$. The
definition of the bracket is the following:  We define the master
symmetry $Y_{-1}$ to be

\bd Y_{-1}=\sum_{i=1}^n  f_i { \p \over \p a_i}  \ , \ed where the
$f_i$ are determined recursively as follows,

\bd f_1=-1, \quad f_{2i}=\frac{a_{2i}} {a_{2i-1}} f_{2i-1} , \quad
f_{2i-1}=-f_{2i-2} -1  . \ed Taking the Lie derivative of $v_2$
in the direction of $Y_{-1}$ we obtain $v_1$, a Poisson bracket
that is homogeneous of degree 1. For $n=5$, $v_1$ takes the
form:

\begin{equation}
\begin{array}{lll}
 \{a_1, a_2\} = a_2 \qquad \{a_1, a_3\} =-a_2  \quad \qquad \{a_1
, a_4\}= { a_2 a_4 \over a_3 } \qquad \{a_1 , a_5 \} =-{ a_2 a_4
\over a_3 } \\
  \{a_2, a_3 \}=a_2 \qquad \{a_2, a_4\} = -{ a_2 a_4 \over a_3
}  \hskip .07cm \qquad \{ a_2, a_5 \}={a_2 a_4 \over a_3 } \label{bracketpi1} \\
 \{a_3, a_4\} = a_4 \qquad \{a_3, a_5 \}=-a_4   \quad \qquad
\{a_4, a_5 \}=a_4.
\end{array}
\end{equation}

Note that the KM Hamiltonian can be expressed as $v_2 dI_1=v_1 dI_2$ and this is a bi-hamiltonian formulation of the KM flow.
\noindent In this paper we rediscover these  brackets using a
recursion operator.

The higher Poisson brackets are  constructed in \cite{damianou2} using a sequence of master symmetries $Y_i$.
For example,  the bracket $v_2$ is
obtained from $v_1$ by taking the Lie derivative in the
direction of the first master symmetry  $Y_1$. Similarly, the Lie derivative of $v_2$ in
the direction of $Y_1$ gives $v_3$.

The brackets $v_1$, $v_2$ and $v_3$ are just the beginning
of an infinite hierarchy which we will derive again in this paper using a different method.

\section{Toda Lattice}
In this paper we will deal with the finite, non-periodic version of the Toda lattice. However, there are also other interesting variations worth mentioning.
There is a  generalization   due to    Deift, Li, Nanda and
 Tomei \cite{deift} who showed that the system remains integrable when $L$ is replaced
 by a full (generic) symmetric $n \times n$ matrix.  The functions  $H_i= { 1 \over i} \, {\rm Tr} \ L^i$  are still in involution but they are  not enough to
ensure integrability. In other words, the existence of a
 Lax pair does not guarantee integrability. There are, however, additional integrals which are rational
functions of the entries of $L$.  The method used to obtain these additional integrals is called chopping and was used originally in
\cite{deift}

The classical Toda lattice was also generalized in another direction.  One can define a Toda type system
for each  simple Lie algebra. The  finite, non--periodic
Toda lattice  corresponds to a root system of type $A_n$. This
generalization is due to Bogoyavlensky \cite{bogo3}. These systems were studied extensively in
 \cite{kostant} where the solution of the system was connected intimately with the representation
theory of simple Lie groups. There are also studies by
 Olshanetsky and Perelomov
\cite{olshanetsky} and Adler, van Moerbeke \cite{avm}.
The description of the systems, following \cite{avm} and \cite{kostant} is as follows:

Let $\g$ be any semi-simple Lie algebra, equipped with its Killing form $\langle\cdot\,\vert\,\cdot\rangle$. One chooses  a Cartan subalgebra
$\h$ of $\g$, a root system $\Delta = \Delta(\h,\g)$ of $\h$ in $\g$,
 a basis of simple roots $\Pi$, and a set of positive roots
 $\Delta^+$.

The Lax pair ($L(t), B(t)$) in ${\g}$ can be described in terms
of the root system as follows:

\bd
L(t)=\sum_{i=1}^l b_i(t) H_{\alpha_i} + \sum_{i=1}^l a_i(t) (X_{\alpha_i}+X_{-\alpha_i}) \ ,
\ed

\bd
B(t)=\sum_{i=1}^l a_i(t) (X_{\alpha_i}-X_{-\alpha_i})  \
\ed

where $H_{\alpha_i}$ is an element of $\h$ and $X_{\alpha_i}$ is a  root vector corresponding to
the simple root $\alpha_i$.

We consider now the classical Toda lattice but we will use a Lax pair slightly different than (\ref{jacobi}).

 Let $D$ be the diagonal matrix with entries $d_i$ where $d_1=1$ and $d_i=a_1 a_2 \dots a_{i-1} \ \ \ i=2,3, \dots N$. In \cite{kostant} Kostant conjugates the matrix $L$ in (\ref{jacobi})  by the matrix $D$ and the resulting matrix $DL D^{-1}$ has the form

\be \label{kostant}
 \begin{pmatrix} b_1 &  1 & 0 & \cdots & \cdots & 0 \cr
                   a_1 & b_2 & 1 & \ddots &    & \vdots \cr
                   0 & a_2 & b_3 & \ddots &  &  \vdots \cr
                   \vdots & \ddots & \ddots & \ddots & & 0 \cr
                   \vdots & & & \ddots & \ddots & 1 \cr
                   0 & \cdots &  \cdots & 0 & a_{N-1} & b_N   \end{pmatrix}  \ .
\ee

Denoting this matrix again by $L$  the equations take the Lax  form

\bd
\dot L(t)=[ L(t), P \, L(t)]
\ed
where $P$ is the projection onto the strictly lower triangular
part of $L(t)$. The decomposition here is strictly lower plus upper triangular.
This form is convenient in applying Lie theoretic techniques to
describe the system.
Note that the diagonal elements correspond to the Cartan subalgebra while the subdiagonal elements correspond to the set $\Pi$ of simple roots.
The equations of motion for the Toda lattice (in Kostant form) are:

\be
\begin{array}{lcl}
 \dot a _i& = & a_i \,  (b_{i+1} -b_i )    \\
   \dot b _i &= &  \,  a_i - a_{i-1}  \ .  \label{f22}
\end{array}
\ee

The functions $ H_i={1 \over i} {\rm Tr} \, L^i$ are independent invariants in involution.

There exists a Lie-Poisson bracket given by the formula
\be
\begin{array}{lcl}
\{a_i, b_i \}& =&-a_i  \\
\{a_i, b_{i+1} \} &=& a_i  \label{f4}   \ ;
\end{array}
\ee
all other brackets are zero.
$H_1=b_1+b_2 + \dots +b_N$ is the only Casimir.   The Hamiltonian in this bracket is
 $H_2 = { 1 \over 2}\  { \rm tr}\  L^2$.  We also have involution of invariants,  $ \{  H_i, H_j \}=0$.
 The Lie algebraic interpretation of this bracket can be found in \cite{kostant}.
  We   denote this bracket by $\pi_1$.

The quadratic Toda bracket appears in conjunction with isospectral deformations of
Jacobi matrices.
 It is a Poisson bracket in which the Hamiltonian vector field generated by $H_1$ is the
same as the Hamiltonian vector field generated by $H_2$ with respect to the $\pi_1$ bracket.  The defining relations are

\be
\begin{array}{lcl}
\{a_i, a_{i+1} \}&=& a_i a_{i+1} \\
\{a_i, b_i \} &=& -a_i b_i                    \\
\{a_i, b_{i+1} \}&=& a_i b_{i+1}    \\
\{b_i, b_{i+1} \}&=&  a_i  \ ;  \label{pi2}
\end{array}
\ee
all other brackets are zero.   The bracket $\pi_2$ is easily defined by taking  the Lie derivative of $\pi_1$  in the direction of  suitable master symmetry  $X_1$, see \cite{damianou1} for details. This bracket has ${\rm det} \, L$ as Casimir and $H_1 ={\rm tr}\, L$
 is the Hamiltonian. The eigenvalues of $L$ are still in involution.
Furthermore, $\pi_2$ is compatible with $\pi_1$.
 We also have
\bd
\pi_2 d H_l = \pi_1   d H_{l+1}  \ . 
\ed

\smallskip
Finally,  we remark that taking the derivative of $\pi_2$  in the direction of $X_1$
yields another Poisson bracket, $\pi_3$, which is cubic in the coordinates.  The defining
relations for $\pi_3$ are
\be
\begin{array}{lcl}
\{a_i, a_{i+1} \}&=&2 a_i a_{i+1} b_{i+1}   \\
\{a_i, b_i \}&=&     -a_i b_i^2-a_i^2 \\
\{a_i, b_{i+1} &=&     a_i b_{i+1}^2 +a_i^2 \\
\{a_i, b_{i+2} \}&=&     a_i a_{i+1} \\
\{a_{i+1}, b_i \}&=&     -a_i a_{i+1} \\
\{b_i, b_{i+1} \} &=&      a_i \, (b_i+b_{i+1}) \ ; \label{pi3}
\end{array}
\ee
all other brackets are zero. The bracket $\pi_3$ is compatible with both $\pi_1$ and
$\pi_2$ and the eigenvalues of $L$ are still in involution.  The Casimir for this
bracket is $ {\rm tr}\, L^{-1}$.

In fact there is an infinite hierarchy of Poisson tensors $\pi_i$,  master symmetries $X_i$ and invariants $H_i$ and they obey some deformation relations.
We quote the results from    refs. \cite{damianou1}, \cite{damianou3}.

\begin{theorem} \label{toda-ab}
\hfill \\
\smallskip
\noindent
{\it i) } $\pi_j$, \  $j\ge 1$ are all Poisson.

\smallskip
\noindent
{\it ii) } The functions $H_i$,        $i\ge 1$ are in involution
 with respect to all of the $\pi_j$.

 \smallskip
 \noindent
 {\it iii)}  $X_i (H_j) =(i+j) H_{i+j} $ ,  $i\ge -1$, $j\ge 1$.

 \smallskip
 \noindent
{\it iv)} $L_{X_i} \pi_j =(j-i-2) \pi_{i+j} $,  $i\ge -1$, $j\ge 1$.

\smallskip
\noindent
{\it v)} $ [X_i, \ X_j]=(j-i)X_{i+j}$, $i\ge 0$, $j\ge 0$.

\smallskip
\noindent
{\it vi)} $\pi_j\  d H_i =\pi_{j-1}\   d  H_{i+1} $.
\end{theorem}

\section{Toda lattice in $(q,p)$ coordinates}

Let ${J}_1$ be the symplectic bracket  with Poisson  matrix

\bd
J_1 = \begin{pmatrix}  0 &  I \cr
                      -I &   0  \end{pmatrix} \ ,
\ed

where $I$ is the $N \times N$ identity matrix. The
bracket $J_1$ is mapped precisely onto the bracket $\pi_1$ under the Flaschka transformation $F$:
\be \label{ft} a_i=e^{q_i- q_{i+1} },  \qquad b_i =-p_i \ . \ee

We define the following tensor,  ${J}_2$, due to Das and Okubo \cite{das}:
\be \label{j2-bracket}
J_2 = \begin{pmatrix} A &  B \cr
                      -B &  C \end{pmatrix}    \ ,
\ee
where $A$ is  the skew-symmetric matrix defined by $a_{ij}=1=-a_{ji}$ for $i<j$,  \,  $B$ is the diagonal matrix
$(-p_1, -p_2, \dots, -p_N)$ and $C$ is the skew-symmetric matrix whose non-zero terms are
$c_{i,i+1}=-c_{i+1,i}=e^{q_i-q_{i+1}}$ for $i=1,2, \dots, N-1$.
The
bracket $J_2$ is mapped precisely onto the bracket $\pi_2$ under the Flaschka transformation (\ref{ft}).

It is easy to see that we have a bi-Hamiltonian pair. We define
\bd
h_1=-(p_1+p_2+\dots +p_N) \ ,
\ed
and $h_2$ to be the Hamiltonian:
\bd
h_2=\sum_{i=1}^N \,  { 1 \over 2} \, p_i^2 +
\sum _{i=1}^{N-1} \,  e^{ q_i-q_{i+1}}  \ .
\ed

Then we obtain the bi-hamiltonian pair

\bd
J_1 d  h_2= J_2 d  h_1 \ .
\ed

We define the recursion operator as follows:

\bd
{\mc R}=J_2 J_1^{-1} \ .
\ed

The matrix form of ${\mc R}$ is quite simple:

\be
{\mc R} ={ 1 \over 2} \begin{pmatrix} B &-A \cr
                       C& B  \end{pmatrix}  \ .     \label{c1}
\ee

Using the recursion operator  we obtain the  higher order Poisson tensors
\bd
J_i = {\mc R}^{i-1} J_1  \ \ \ i=2, 3, \dots \ .
\ed

Following \cite{fern1} we   define the conformal symmetry
\bd
Z_0=\sum_{i=1}^N  (N-2i+1) {\partial \over \partial q_i} +\sum_{i=1}^N p_i {\partial \over \partial p_i} \ .
\ed

It is straightforward to verify that
\bd
{\mc L}_{Z_0} J_1=- J_1 \ ,
\ed

\bd
{\mc L}_{Z_0} J_2=0  \ .
\ed

In addition,
\bd
Z_0(h_1)=h_1
\ed

\bd
Z_0(h_2)=2h_2  \ .
\ed

Consequently, $Z_0$ is a conformal symmetry for $J_1$, $J_2$ and $h_1$. The constants appearing in Theorem \ref{oevel} are
$\lambda=-1$, $\mu=0$ and $\nu=1$. According to Oevel's Theorem
 we end up with the following deformation relations:

\bd
[Z_i, h_j]= (i+j)h_{i+j}
\ed

\bd
L_{Z_i}  J_j = (j-i-2) J_{i+j}
\ed

\bd
 [ Z_i, Z_j ]  = (j-i) Z_{i+j}  \ .
\ed

Switching to Flaschka coordinates, we obtain  relations  iii)- v) of Theorem \ref{toda-ab}.

\section{From Toda to Volterra }

We consider $T_{(a,b)}$ the phase space of the Toda lattice in Flaschka coordinates and the space $V_a$ of KM-system in $a$ coordinates.
Note that $V_a$ is not a
Poisson subspace of $T_{(a,b)}$.  However, as was demonstrated in \cite{damianou5}   $V_a$ is the fixed manifold of
the involution $\phi:T_{(a,b)}\to T_{(a,b)}$ defined by
\[   \phi (a_1,a_2\dots,a_{N-1},b_1,b_2\dots,b_N)\mapsto
(a_1,a_2\dots,a_{N-1} ,-b_1,-b_2\dots,-b_N),\]
and we have the following result:

\begin{theorem}
\label{prop:A:Volterra}
$\phi:T_{(a,b)}\to T_{(a,b)}$ is a Poisson automorphism of $(T_{(a,b)}, \pi_k)$, if
$k$ is even.
\end{theorem}

Therefore, by Corollary \ref{cor:Poisson:involution}, $V_a$
inherits a family of Poisson brackets $v_2, v_3,\dots$. For example,
the quadratic bracket can be computed from formulas
(\ref{pi2})  and is given by
\begin{equation}
\label{eq:quadratic:Volterra}
 \pb {a_i}{a_j}=a_ia_j(\delta_{i,j+1}-\delta_{i+1,j}),
\end{equation}
while the Poisson bracket $v_3$ is found to be given by the
formulas
\begin{align}
\label{eq:cubic:Volterra}
 \pb {a_i}{a_{i+1}}&=a_ia_{i+1}(a_i+a_{i+1}),& (i&=1,\dots,n-1)\\
\pb  {a_i}{a_{i+2}}&=a_ia_{i+1}a_{i+2},& (i&=1,\dots,n-2) \notag
\end{align}

i.e. we obtain $v_2$ and $v_3$ of Section 4. It follows that the restriction of the integrals $H_{2k}$ to $V_a$
gives a set of commuting integrals, with respect to these Poisson
brackets. Also,  the Lax equations (\ref{kostant}) lead
to Lax equations for the corresponding flows, merely by putting all
$b_i$ equal to zero. We  recover the vector field
\begin{equation}\label{KM}
  \dot{a}_i=a_i(a_{i-1}-a_{i+1}),\qquad i=1,\dots,n,
\end{equation}
and  a family of integrable
systems admitting a multiple hamiltonian formulation:
 \bd v_k dI_{l}=v_{k-1} dI_{l+1}, \qquad (k=1,2,\dots) , \ed
where $v_k$ is the restriction of $\pi_{2k}$ to  $V_a$ and $I_l$ is the restriction of $H_{2l}$ to $V_a$.

\section{Symplectic realization }

Let $N=2(n+2)$ and consider  coordinates $(q_1, q_2, \dots, q_N)$ in ${\bf R}^N$.
We define the following transformation $G$  from ${\bf R}^{2n+2}$ to
${\bf R}^{2n+1}$,

\be  \label{b1}
a_i=e^{ q_i -q_{i+1}}  \qquad i=1, 2,  \dots, N-1  \ . \end{equation}

\noindent The Hamiltonian in $q$ coordinates is given by

\be
\label{b2}
i_1=  \sum_{i=1}^{N-1}
e^{q_{i}-q_{i+1}}  \ .
  \ee

 We define the Poisson bracket
 \bd
 \{ q_i, q_j \}=1  \qquad \forall \ \  i < j  \ .
 \ed
 Let us denote this constant Poisson tensor  by $w_2$ (since it is the analogue of $v_2$ in $q$ coordinates).

 Hamilton's equation $w_2 di_1$  become
 \be
 \dot{q}_i =-e^{q_{i-1}-q_i} -e^{q_i -q_{i+1} } \ . \label{q-equations}
 \ee

\noindent It is straightforward to check that Hamilton's equations
 correspond in the $a-$space to the KM-system
(\ref{a1}) via the mapping $G: \R^N  \to \R^{N-1}$

\be G(q_1, \dots, q_N) =\left( e^{q_1- q_2}, \dots, e^{q_{N-1}- q_N} \right) \ . \label{g-map} \ee

 In fact we calculate:
\begin{align*}
\dot{a}_i
&  =e^{q_i -q_{i+1}} (\dot{q}_i -\dot{q}_{i+1} ) \\
&=a_i \left( -e^{q_{i-1} -q_i}- e^{q_i -q_{i+1}} +e^{q_i -q_{i+1} } + e^{q_{i+1}-q_{i+2}} \right) \\
&= a_i(a_{i+1} -a_{i-1}) \ . \end{align*}

The symplectic
bracket $w_2$ in $V_q$   space  corresponds to the quadratic
bracket $v_2$ in $V_a$ space.

We then
define a bracket $w_3$ in ${\bf R}^{N}$ which is mapped to
$v_3$ under the transformation $G$. The  bracket $w_3$ is defined by the formula:

\bd \{ q_i, q_j \}=e^{q_{i-1}-q_i}+(1-\delta_{i+1, j}) e^{q_i- q_{i+1}} + e^{q_{j-1}-q_j}+e^{q_j-q_{j+1}} \ . \ed
Whenever a term is not defined we omit that whole term. This happens only when $i=1$ or $j=N$.  We note that $w_2$ is
compatible with $w_3$.

Define \bd i_0= \sum_{k=1}^N (-1)^{k+1} q_k=q_1-q_2 + q_3-q_4+ \dots +q_{N-1}-q_N  \ . \ed

\begin{rem}
It is not difficult to discover $i_0$. According to \cite{damianou7},  in the presence of an invertible
Nijenhuis tensor a natural choice of functions to form a bi-hamiltonian pair is

$\frac{1}{2}\log(\det\mc{R}) $ and $\frac{1}{2}\tr\mc{R}  .$  It turns out that the determinant of the Toda recursion operator restricted to the Volterra phase space is $e^{2 i_0}$  and the trace equals $2i_1$.

\end{rem}

  Then
\bd \dot{q}_i =\{q_i, i_0 \}_{w_3}  =-e^{q_{i-1}-q_i }-e^{q_i -q_{i+1}} \ . \ed
In other words we have a bi-hamiltonian pair

\bd w_2 di_1= w_3 di_0 \ . \ed

 We define a recursion operator as follows:

\begin{equation*}
{\mathcal R}=w_3 w_2^{-1} .
\end{equation*}
 In $q$ coordinates, the  symbol
$\chi_j$ is a shorthand for  $\chi_{i_j}$. It is generated, as
usual, by

\begin{equation*} \chi_i = {\mathcal R}^{i-1} \chi_1 . \end{equation*}

Note that $i_1$ corresponds under mapping (\ref{g-map}) to a
constant multiple of $I_1=\frac{1}{2} {\rm Tr} \, (L)^2$. In a
similar fashion we obtain the higher order Poisson tensors

\bd w_i = {\mathcal R}^{i-2} w_2 \qquad i=3,4, \dots . \ed

\noindent We finally define the conformal symmetry \bd
X_0=\sum_{i=1}^N (N-i+1)  {\partial \over
\partial q_i}  \ . \ed

\noindent The Poisson tensors $w_2, w_3$ and the functions
$i_0,i_1$ define a bi-Hamiltonian pair.    It is straightforward to
verify that

\begin{equation*} {\mathcal L}_{X_0} w_2=0, \quad
{\mathcal L}_{X_0} w_3= w_3, \quad {\mathcal L}_{X_0} i_1=i_1 .
\end{equation*}
Consequently, $X_0$ is a conformal symmetry for $w_2$, $w_3$ and
$i_1$. The constants appearing in Oevel's Theorem are $\lambda=0$,
$\mu=1$ and $\nu=1$. Therefore, we end up with the following
deformation relations:

\bd [X_k, i_j]= (k+j) i_{k+j} \ed

\bd L_{X_k}  v_j = (j-k-2) v_{k+j} \ed

\bd
 [ X_k, X_j ]  = (j-k) X_{k+j}  \ .
\ed

\noindent Projecting to the $a-$space under mapping
(\ref{g-map}) we obtain the multiple hamiltonian structures of \cite{damianou2}.

\begin{rem}
Note that we may define a Poisson tensor  $w_1$ by the formula $w_1=w_2 w_3^{-1} w_2$ as in \cite{damianou4}.  The projection of $w_1$ under transformation (\ref{g-map}) gives $v_1$ of Section 4.
Perhaps a more appropriate  bottom arrow in the diagram of the introduction is
\bd
 G: (w_1, w_2, i_1, i_2) \to (v_1, v_2, I_1, I_2) \ .
 \ed
\end{rem}

\section{From Toda $(q,p)$ space to Volterra $q$ space}

In this Section we will explain the origin of the symplectic bi-hamiltonian realization of Section 8. The  idea  is to use a Poisson involution in the Toda $(q,p)$ space and to reduce the equations to the Volterra $V_q$ space.
We consider  $T_{ (q,p)}$, the phase space of the Toda lattice in natural coordinates.  In Section 6 we have constructed a bi-hamiltonian system given by the Poisson tensors $J_1$, $J_2$ and the Hamiltonians $h_1$ and $h_2$. We also have a recursion operator ${\mc R}$   which gives rise to a sequence of Poisson tensors $J_i \ \ i=1,2, \dots $.

Define the involution  $\psi : T_{ (q,p)} \to T_{ (q,p)}$ by the formula

\bd \psi (q_1, \dots, q_N, p_1, \dots, p_N)=(q_1, \dots, q_N, -p_1, \dots, -p_N) \  \ . \ed

We have the following result:

\begin{theorem}
\label{prop:B:Volterra}
$\psi: T_{ (q,p)} \to T_{ (q,p)}  $ is a Poisson automorphism of $( T_{ (q,p)}, J_{2k})$,  $k=1,2,\dots $.
\end{theorem}

We will not give the proof of this result since it is entirely analogous to the proof of Theorem \ref{prop:A:Volterra} which is given in \cite{damianou5}.

Therefore using  Corollary \ref{cor:Poisson:involution}, $V_q$ inherits a family of Poisson tensors $w_2, w_3, \dots$. For example, the bracket $w_2$ is clearly the $A$ block of the Poisson matrix $J_2$ in (\ref{j2-bracket}).  On the other hand it is straightforward to compute that in the bracket $J_4$  we have
\bd \{q_i, q_j \}=p_i^2 +p_i p_j +p_j^2 + e^{q_{i-1}-q_i}+(1-\delta_{i+1, j}) e^{q_i- q_{i+1}} + e^{q_{j-1}-q_j}+e^{q_j-q_{j+1}} \ . \ed

Therefore the reduction of $J_4$ to $V_q$ is precisely $w_3$ of Section 8.

\section{Moser's recipe}

Moser in \cite{moser2} describes a relation between the KM system  (\ref{a1}) and the non--periodic Toda lattice. The procedure is
the following: Form $L^2$  which is not anymore a tridiagonal matrix but is similar to one.  Let $\{e_1, e_2, \dots, e_n \}$ be the standard
basis of ${\bf R}^n$, and $E_o= \{ {\rm span}\, e_{2i-1}, \,  i=1,2, \dots \}$, $E_e= \{ {\rm span}\, e_{2i}, \, i=1,2, \dots \}$. Then $L^2$ leaves
$E_o$,  $E_e$ invariant and reduces to each of these spaces to a tridiagonal symmetric Jacobi matrix.
 For example, if we omit  all even columns and all even  rows we
obtain a tridiagonal Jacobi matrix and the entries of this new matrix  define the transformation from the KM--system
to the Toda lattice. We illustrate with a simple example where $n=5$.

We use the symmetric version of the KM system Lax pair (due to Moser)  given by
\be
L=\begin{pmatrix} 0 & a_1& 0 & 0& 0\cr
                 a_1&0& a_2& 0& 0\cr
                 0 & a_2 &0 & a_3 & 0  \cr
                 0 & 0 & a_3& 0 &a_4 \cr
                  0 &0&0&a_4& 0
                  \end{pmatrix}  \ .
\ee

It is simple to calculate that  $L^2$ is the matrix

\be
\begin{pmatrix} a_1^2 & 0& a_1 a_2 & 0& 0 \cr
                 0& a_1^2+a_2^2& 0& a_2 a_3&  0\cr
                 a_1 a_2 & 0 & a_2^2+a_3^2 & 0 & a_3 a_4\cr
                 0 & a_2 a_3 & 0& a_3^2+a_4^2 &0\cr
                  0 &0& a_3 a_4&0&a_4^2
                  \end{pmatrix} \ .
\ee
Omitting even  columns and even rows of $L^2$ we obtain the matrix
\be
\begin{pmatrix}  a_1^2 & a_1 a_2 & 0 \cr
                   a_1 a_2 &a_2^2+a_3^2 & a_3 a_4 \cr
                    0 & a_3 a_4 & a_4^2
                  \end{pmatrix} \ .
\ee
This is a tridiagonal Jacobi matrix. It is natural to define   new variables $A_1=a_1 a_2$, $A_2=a_3 a_4$, $B_1=a_1^2$, $B_2=a_2^2+a_3^2$, $B_3=a_4^2$. The new
 variables $A_1,A_2, B_1,B_2, B_3$ satisfy the Toda lattice  equations (\ref{toda-tri}).

This procedure shows that the KM-system  and the Toda lattice are closely related: The explicit  transformation
 which is  due to H\'enon
maps one system to the other. The mapping in the general case  is given by
\be
A_i=-{ 1 \over 2} \sqrt {a_{2i} a_{2i-1}} \  , \qquad  B_i= { 1 \over 2}\left( a_{2i-1}+a_{2i-2} \right)  \label{a25} \ .
\ee
The equations satisfied by the new variables $A_i$, $B_i$ are given by:

\begin{displaymath}
\begin{array}{lcl}
 \dot A _i& = & A_i \,  (B_{i+1} -B_i )    \\
   \dot B _i &= & 2 \, ( A_i^2 - A_{i-1}^2 ) \ .
\end{array}
\end{displaymath}
These are precisely the Toda equations in Flaschka's form \cite{flaschka1}.

This idea of Moser was applied with success to establish transformations   from the generalized Volterra  lattices of Bogoyavlensky \cite{bog1, bog2}  to generalized Toda systems.
The relation between the Volterra systems of type $B_n$ and $C_n$ and the corresponding Toda systems is in \cite{damianou5}.  The similar construction of  the Volterra lattice   of type $D_n$ and the generalized
Toda lattice of type $D_n$   is in \cite{damianou6}.  It turns out  that the Volterra $D_n$ system
corresponds not to the Toda $D_n$ system but to a special case of the Sklyanin lattice \cite{sklyanin}.

\end{document}